\documentclass[review]{elsarticle}

\usepackage{hyperref}

\usepackage{verbatim}
\usepackage[utf8]{inputenc}
\usepackage[english]{babel}
\usepackage{amssymb}
\usepackage{amsthm}
\usepackage{amsmath,amssymb,amsfonts,graphicx,float,multicol,fleqn}
\usepackage{subcaption}
\usepackage{multirow}

\usepackage{algorithm}
\usepackage[noend]{algpseudocode}
\newcommand{\etal}{{\em et~al.\/}}
\newcommand{\IR}{\ensuremath{\mathbb{R}}} 
\newcommand{\expnumber}[2]{{#1}\mathrm{e}{#2}}

\usepackage[letterpaper, margin=1in]{geometry}
\usepackage[dvipsnames]{xcolor}

\newcommand{\RN}[1]{\textup{\uppercase\expandafter{\romannumeral#1}}}

\renewcommand{\paragraph}[1]{\bigskip\noindent\textbf{#1} }

\newcommand{\REM}[1]{}

\newcommand{\CostFun}{\ensuremath{\mathop{\mathrm{CostFun}}}}

\begin{document}

\begin{frontmatter}

\title{Legendre Deep Neural Network (LDNN) and its application for approximation of nonlinear Volterra–Fredholm–Hammerstein integral equations}

\author[1]{Zeinab Hajimohammadi}
\ead{Z\_Hajimohammadi@sbu.ac.ir}
\author[1,3,4]{Kourosh Parand}
\ead{k\_parand@sbu.ac.ir}
\author[4]{Ali Ghodsi}
\ead{aghodsib@uwaterloo.ca}

\address[1]{Department of Computer and Data Sciences, Faculty of Mathematical Sciences, Shahid Beheshti University, Tehran, Iran}
\address[3]{Institute for Cognitive and Brain Sciences, Shahid Beheshti University, Tehran, Iran}
\address[4]{Department of Statistics and Actuarial Science, University of Waterloo, Waterloo, Canada}

\begin{abstract}
Various phenomena in biology, physics, and engineering are modeled by differential equations. These differential equations including partial differential equations and ordinary differential equations can be converted and represented as integral equations. In particular, Volterra–Fredholm–Hammerstein integral equations are the main type of these integral equations and researchers are interested in investigating and solving these equations.  In this paper, we propose Legendre Deep Neural Network (LDNN) for solving nonlinear Volterra–Fredholm–Hammerstein integral equations~(V-F-H-IEs). LDNN utilizes Legendre orthogonal polynomials as activation functions of the Deep structure. We present how LDNN can be used to solve nonlinear V-F-H-IEs. We show using the Gaussian quadrature collocation method in combination with LDNN results in a novel numerical solution for nonlinear V-F-H-IEs. Several examples are given to verify the performance and accuracy of LDNN.
\end{abstract}


\end{frontmatter}

\section{Introduction}
Various phenomena in biology, physics, finance, neuroscience and engineering are modeled by differential equations \cite{gandzha2019high, atangana2020fractional, courant2008methods,heydari2019chebyshev, davis1961introduction}. Integral models and  application of them have been observed in many fields of sciences, mathematics,  engineering and technology and research works in this field have expanded greatly \cite{bazm2020bernoulli,beiglo2020numerical,babolian2007chebyshev,amiri2020spectral,xie2020applying, sarwarexistence,mohammad2020collocation,maleknejad2020convergence}.\\The various numerical methods are applied for solving differential and integral equations. Homotopy analysis method (HAM) \cite{liao2012homotopy} and variational iteration method (VIM) \cite{he2007variational} are known as analytical/semi-analytical methods. Usually, Nystrom method \cite{lardy1981variation}, spectral methods \cite{mohammad2020collocation,canuto2012spectral,parand2016solving,marzban2011composite,parand2012numerical}, Runge-Kutta methods \cite{hairer2006numerical}, the finite difference methods (FDM) \cite{smith1985numerical} and the finite element methods (FEM) \cite{johnson2012numerical} are considered as the popular numerical methods. When the complexity of the model does not allow us to obtain the solution explicitly, numerical methods are a proper selection for finding the approximate solution for the models.
Deep neural networks are a main and beneficial part of machine learning family which are applied in various areas including engineering, speech processing, computer vision, natural language processing and image processing and etc. \cite{li2019deep, peiris2021generalised,lecun2015deep,dang2020multilayer,krizhevsky2012imagenet,gorgel2019face}.
Also, the approximation of the functions is a significant branch in scientific computational  and  achieving success in this area  is considered  by some research \cite{mohamadipanah2021deep, tang2019chebnet, hanin2019universal}.
Solving differential equations is the other main branch of scientific computational which neural networks and deep learning have been shown success in this area.  \cite{lample2019deep,berg2018unified,zhang2019bilinear,raissi2019physics,kuhnel2019differential}. 
 In recent years, several researchers studied the solving differential equations via deep learning or neural networks. differential equations consists of ordinary differential equations, partial differential equations and integral equations. \cite{sirignano2018dgm,lu2019deepxde,meng2020ppinn,kuhnel2019differential}.
 
Recently, some of the machine learning methods are applied for solving differential equations. Chakraverty and Mall \cite{chakraverty2017artificial}
introduced orthogonal neural networks which used orthogonal polynomials in the structure of the network.  Raja \etal \cite{raja2019numerical} applied meta-heuristic optimization algorithm to neural network for obtaining the solution of differential equations. Moreover, other methods of machine learning such as support vector machine \cite{vapnik2013nature} are used to approximate the solution of the models. Least squares support vector machines are considered in these researches \cite{hajimohammadi2020new, mehrkanoon2015learning}.
Baker \etal \cite{baker2019workshop} selected deep neural networks for solving the differential equations. Pang \etal \cite{pang2019fpinns} introduced a new network to find the solution of the different equations. Han \etal \cite{han2018solving} solved high-dimensional problems via deep networks. Also, Long \etal \cite{long2018pde} and Raissi \etal \cite{raissi2019physics} introduced a group of the equations which solved by deep learning.
Furthermore,  He \etal \cite{he2018relu} and Molina \etal \cite{molina2019pade} investigated the effect of the activation function on networks.\\
In this paper, we concern nonlinear Volterra–Fredholm–Hammerstein integral equations~(V-F-H-IEs) and try to obtain the solution of them via deep neural network. We present a new numerical approach of machine learning which is a combination of deep neural network and Legendre collocation method. This approach is useful for solving the differential equations and we applied it for solving nonlinear V-F-H-IEs. We used Legendre collocation method to our network for perfect the numerical computations and enhancement the performance the network.

\section{Legendre Deep Neural Network (LDNN)}
The main purpose of introducing LDNN is to apply it for solving differential models. Indeed, this purpose is to expand the utilization of deep learning networks in the field of scientific computing, especially the solution of differential equations. Moreover, this network has the advantages of solving equations by deep learning as well as numerical methods such as collocation method used to achieve better solution to the equations. LDNN presents a combination of a deep neural network and Legendre collocation method. In fact, our network consists of two networks  which have connected consecutive to each other. The first network is a feed forward neural network which has an orthogonal Legendre layer. The second network includes operation nodes to create the desired computational model. In recent decades, numerical methods especially collocation method are popular methods for solving differential equations. In the collocation method,  first an approximation of the solution is expanded by using the sum of the basic functions.  The basic functions consists of the orthogonal polynomials such as Legendre polynomials.Then this approximation is placed in the differential equation. By considering the appropriate set of candidate points, an attempt is made to obtain the unknown coefficients of the basic functions so that  the solution  satisfies  the equation in a set of candidate points. The first network is applied to creat the approximation of the solution. This approximation can be known as the scattered data interpolation problem. The second network is used to obtain the desired equation so that the solution satisfies it. The structure of LDNN is described in detail at the following rest.\\
Consider that the first network has a $\mathcal{M}$-layer which defined as follows:
\begin{eqnarray}
&& \mathcal{H}_0 = x, \quad x \in \IR^d, \nonumber\\
&& \mathcal{H}_1 = L (W^{(1)}  \mathcal{H}_0 + b^{(1)}), \nonumber\\
&& \mathcal{H}_i = f (W^{(i)}  \mathcal{H}_{i-1} + b^{(i)}),\quad 2 \leq i \leq \mathcal{M}-1, \nonumber\\
&& \mathcal{H}_{\mathcal{M}} = W^{(\mathcal{M})}  \mathcal{H}_{\mathcal{M}-1}+b^{(\mathcal{M})}. \nonumber
\end{eqnarray}
where $\mathcal{H}_0$ is the input layer with $d$ dimension. $\mathcal{H}_i,\; 1\leq i \leq \mathcal{M}-1$ are hidden layers, $L=[L_0, L_1,...L_n]^T$ which $L_i$ are \textit{i}-th degrees of Legendre orthogonal polynomials, $\mathcal{H}_1$ is an orthogonal layer, $f$ is the hyperbolic tangent activation function or other commonly used activation functions.
$W^{(i)}, \; i= 1, \cdots , \mathcal{M} $ are the weight parameters and $b^{(i)}, \; 1 \leq i \leq \mathcal{M}$ are the  bias parameters.
$\mathcal{H}_{\mathcal{M}}$ is the output layer. It is notable that the second network is applied to obtain the desired differential model. This aim is possible by using operation nodes including integrals, derivatives, and etc. These nodes are applied to the output of the first network. Moreover, automatic differentiation~(AD)~\cite{baydin2017automatic} and Legendre Gaussian integration~\cite{shen2011spectral} have been used in network computing to obtain more accurate and fast calculations. How to train the network and set the parameters are also important points. Supervised learning method is used to train network. The cost function for setting parameters is defined as follows:
\begin{equation}
    \CostFun =  \min( y_t- y_p) + \min(R_m).
\end{equation}
where $ y_t$  is an exact value of the model and $ y_p$ is a predicted value of the LDNN. The definition of $R_m$ is explained in section~\ref{V_F_H_IEs}.The minimization of $\CostFun$ is obtained by performing Adam algorithm~\cite{kingma2015adam} and the L-BFGS method~\cite{liu1989limited} on mean squared errors of training data set.
\subsection{Legendre polynomials}
Legendre polynomials~\cite{shen2011spectral} are a main series of orthogonal polynomials which denoted by $L_n(\eta)$, are defined as:
\begin{equation}
    L_n(\eta) = \frac{1}{2^n} \sum_{\ell=0}^{[\frac{n}{2}]} (-1)^{\ell}\frac{(2n-2\ell)!}{2^n\ell!(n-\ell)!(n-2\ell)!}\eta^{n-2\ell}
\end{equation}
Legendre polynomials are defined in $[-1,1]$ domain and have the recurrence formula in the following form:
\begin{eqnarray}
&& (n+1)L_{n+1}(\eta) = (2n+1) \eta L_n(\eta) - nL_{n-1}(\eta), \quad n \geq 1,\nonumber\\
&& L_0(\eta)=1, \quad\quad L_1(\eta)=\eta.
\end{eqnarray}
Orthogonality relation for these polynomials is as follows:
\begin{equation}
    \int_{-1}^{1} L_n(\eta) L_m(\eta)\mathrm{d}\eta = \gamma \delta_{n,m},
\end{equation}
where $\delta_{n,m}$ is a delta Kronecker function and $\gamma=\frac{2}{2n+1}$.\\
The weight function of them is $\mathcal{W}(\eta)=1$.
Some following useful properties of Legendre polynomials are defined:
\begin{eqnarray}
&& L_n(-\eta)=(-1)^nL_n(\eta),\\
&& |L_n(\eta)| \leq 1, \quad \forall \eta \in [-1,1],\;n\geq 0,\\
&& L_n(\pm 1)=(\pm 1)^n,\\
&& (2n+1)L_n(\eta)=L'_{n+1}(\eta)- L'_{n-1}(\eta),\quad n\geq1.
\end{eqnarray}
\section{Nonlinear Volterra–Fredholm–Hammerstein integral equations and LDNN}
\label{V_F_H_IEs}
\begin{figure}[h]
\begin{center}
\includegraphics[width=1.0\textwidth]{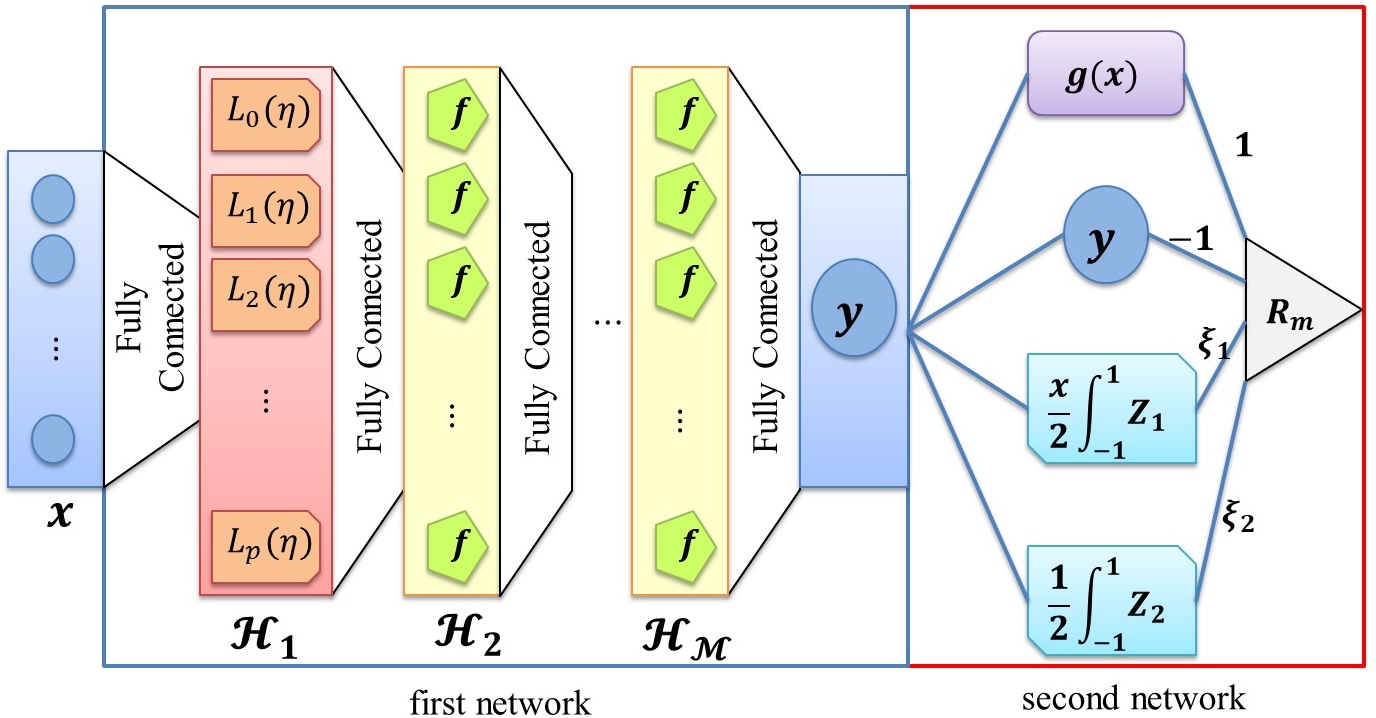}
\end{center}
\caption{The architecture of LDNN for solving nonlinear V-F-H-IEs. The first network approximates the solution of IE $y(x)$. This network has $\mathcal{M}$-layer and feed forward neural network is the structure of it. $\mathcal{H}_1$ is introduced as a orthogonal layer which consists of $p$ neurons with $\{L_i\}_{i=0}^{p}$ (Legendre polynomials) as activation functions. Other layers have $f$, hyperbolic tangent as activation functions. The second network with the nodes makes the desired model and the output of it, is $R_m$ (consider Eq. (\ref{R_m})). The outputs of LDNN are $y(x)$ and $R_m$.}
\label{archi}
\end{figure}

The general form of nonlinear Volterra–Fredholm–Hammerstein integral equations~(V-F-H-IEs) is as follows:
\begin{equation}
\label{VFHIEEQ}
    y(x)=g(x) + \xi_1 \int_{0}^{x} K_1(x,s) \varphi_1(s,y(s)) \mathrm{d}s + \xi_2 \int_{0}^{1} K_2(x,s) \varphi_2(s,y(s)) \mathrm{d}s, \quad x \in [0,1].
\end{equation}
where $\xi_1$, $\xi_2$ are fixed, $g(x)$, $K_1(x,s)$ and $K_2(x,s)$ are given functions and $\varphi_1(s,y(s))$, $\varphi_2(s,y(s))$ are nonlinear functions. The aim is to find the proper $y(x)$. In order to use the LDNN, reformulated Eq. (\ref{VFHIEEQ}) in the following form:
\begin{equation}
\label{res}
    R_m = -y(x) + g(x) + \xi_1 \int_{0}^{x} K_1(x,s) \varphi_1(s,y(s)) \mathrm{d}s + \xi_2 \int_{0}^{1} K_2(x,s) \varphi_2(s,y(s)) \mathrm{d}s, \quad x \in [0,1].
\end{equation}
$y(x)$ is approximated by the first network of the LDNN.
\begin{equation}
    y(x) \approx \mathcal{H}_{\mathcal{M}} .
\end{equation}
Furthermore, we applied Legendre–Gauss integration formula~\citep{shen2011spectral}:
\begin{equation}
    \int_{-1}^{1} h(X)\mathrm{d}X = \sum_{j=0}^{N} \omega_j h(X_j)
\end{equation}
where $\{X_j\}_{j=0}^N$ are the roots of $L_{n+1}$ and $\{\omega_j\}_{j=0}^N = \frac{2}{(1-X^2_j)(L'_{n+1}(X_j))^2}$. Here, we should transfer the  $[0,x]$ and $[0,1]$ domains into the $[-1,1]$ domain. It is possible by using the following transformation:
\begin{equation}
    t_1 = \frac{2}{x}s-1, \quad t_2 = 2s-1. \nonumber
\end{equation}
Consider
\begin{eqnarray}
&& Z_1(x,s)=K_1(x,s) \varphi_1(s,y(s)),\nonumber\\
&& Z_2(x,s)=K_2(x,s) \varphi_2(s,y(s)).\nonumber
\end{eqnarray}
we have
\begin{equation}
    R_m = -y(x) + g(x) + \xi_1 \frac{x}{2} \int_{-1}^{1} Z_1(x,\frac{x}{2}(t_1+1))\mathrm{d}t_1 + \frac{\xi_2}{2} \int_{-1}^{1} Z_2(x,\frac{x}{2}(t_2+1)) \mathrm{d}t_2.
\end{equation}
by using Legendre–Gauss integration formula, the below form is concluded:
\begin{equation}
\label{R_m}
    R_m = -y(x) + g(x) + \xi_1 \frac{x}{2} \sum_{j=0}^{N_1} \omega_{1j} Z_1(x,\frac{x}{2}(t_{1j}+1)) + \frac{\xi_2}{2} \sum_{j=0}^{N_2} \omega_{2j}  Z_2(x,\frac{x}{2}(t_{2j}+1)).
\end{equation}
The second network of LDNN and its nodes makes $R_m$. The architecture of LDNN for solving nonlinear V-F-H-IEs is represented in Figure~\ref{archi}. 

\section{Numerical results}
In order to present the accuracy and performance of the LDNN for solving nonlinear V-F-H-IEs and justify the efficiency of the proposed method, several examples are given. The convergence behavior of the LDNN is reported by using the following parameters:\\
The exact value $y_t$, the predicted value $y_p$ and the absolute error~(Error) in some points of test data are reported in various tables. The number of the train data $m_1$, the number of Legendre quadrature points $(N_1,N_2)$, the number of the test data $m_2$, the structure of network $\mathcal{M}$-layers, $L_{2}^{train}$ and $L_{2}^{test}$ are shown in Table~\ref{param-table}. $L_{2}^{train}$ and $L_{2}^{test}$ are calculated as follows:
\begin{eqnarray}
&& L_{2}^{train}= || y_t - y_p||_2 = [\sum_{j=1}^{m_{tr}}(y_t(x_j)-y_p(x_j))^2]^{\frac{1}{2}}, \nonumber\\
&& L_{2}^{test}= || y_t - y_p||_2 = [\sum_{j=1}^{m_{te}}(y_t(x_j)-y_p(x_j))^2]^{\frac{1}{2}},
\end{eqnarray}
The Tensorflow package of Python version 3.7.0. is applied for writing the code of all experiments. Adam algorithm is stoped when the number of iteration is up to $5000$ and L-BFGS method is stoped when it converges. The figures are obtained on the test data set.For comparison the results of the LDNN with other well-known methods, the Adomian decomposition method~(ADM) \cite{wazwaz2011linear} and a simple feed forward neural network~(FNN) which have supervised training and similar number of parameters and standard loss functions without the IE residual, are considered. The numerical results are reported in the tables.
\begin{table}[t]
\caption{The LDNN parameters  for all the experiments. The structure of $\mathcal{M}$-Layers indicates by $[d,NL^{(1)},NL^{(2)},\cdots,NL^{(\mathcal{M}-1)},1]$. This network has $d$ dimension in input layer, $\mathcal{M}-1$ hidden layers with $NL^{(\ell)},\;2\leq \ell \leq \mathcal{M}-1$, neurons in each layer and one output which approximates the $y(x)$. All the experiments have $4$ hidden layers.}
\label{param-table}
\begin{center}
\begin{tabular}{ccccccc}
\bf Experiment &$\mathcal{M}$-Layers&$m_1$&$(N_1,N_2)$&$L_{2}^{train}$&$m_2$&$L_{2}^{test}$\\
\\ \hline \\
Experiment 1&$[1,10,30,20,10,1]$&$500$&$(50,-)$&$\expnumber{3.937867}{-09}$&$100$&$\expnumber{4.015095}{-09}$\\
Experiment 2&$[1,10,30,20,10,1]$&$500$&$(50,50)$&$\expnumber{7.156029}{-09}$&$100$&$\expnumber{7.537263}{-09}$\\
Experiment 3&$[1,10,30,20,10,1]$&$500$&$(50,50)$&$\expnumber{1.347132}{-09}$&$100$&$\expnumber{1.659349}{-08}$\\
Experiment 4&$[1,10,30,20,10,1]$&$500$&$(50,50)$&$\expnumber{9.182442}{-09}$&$100$&$\expnumber{1.107755}{-09}$\\
\end{tabular}
\end{center}
\end{table}
\subsection{Experiment~1}
\begin{figure}[h]
\begin{center}
\begin{subfigure}[b]{0.3\textwidth}
         \begin{center}
         \includegraphics[width=\textwidth]{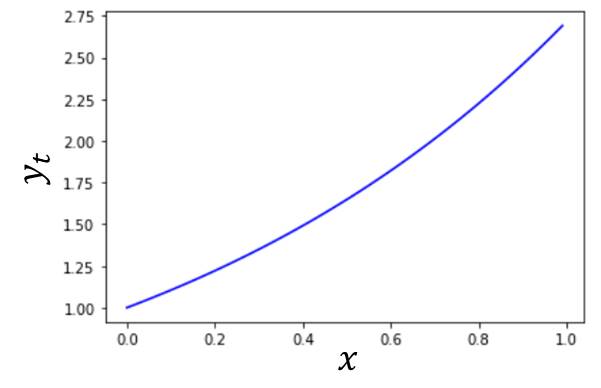}
         \label{fig:1-1}
         \end{center}
     \end{subfigure}
     \begin{subfigure}[b]{0.3\textwidth}
         \begin{center}
         \includegraphics[width=\textwidth]{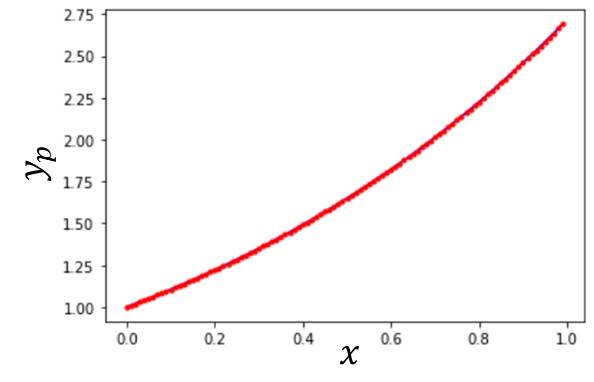}
         \label{fig:1-2}
         \end{center}
     \end{subfigure}
     \begin{subfigure}[b]{0.3\textwidth}
         \begin{center}
         \includegraphics[width=\textwidth]{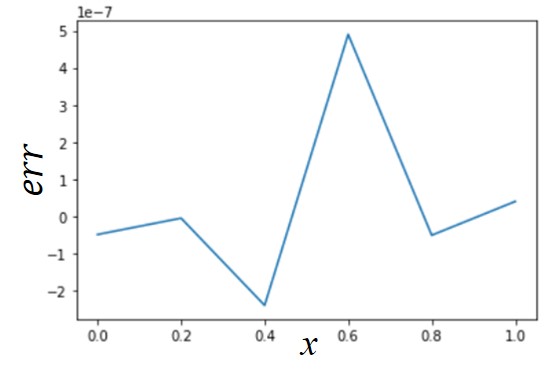}
         \label{fig:1-3}
         \end{center}
     \end{subfigure}
\end{center}
\caption{Results of Experiment~1. Exact solution $y_t(x)=\mathrm{e}^{x}$,  predicted solution $y_p(x)$ by LDNN and absolute error $err=|y_t(x)-y_p(x)|$.}
\label{fig-Ex1}
\end{figure}

Suppose that we have the following model \cite{yousefi2005legendre}:
\begin{equation}
\label{Ex1}
    y(x)=\mathrm{e}^{x} -\frac{1}{3}\mathrm{e}^{3x} + \frac{1}{3} +\int_{0}^{x} y^3(s) \mathrm{d}s, \quad x \in [0,1].
\end{equation}
It has the exact solution $y(x)=\mathrm{e}^{x}$. Table~\ref{Ex1-table} represents the exact value, the predicted value and the absolute error~(Error) in several test points on $[0,1]$ domain. $50$ points of shifted Legendre quadrature points are applied for training LDNN. The number of train data set is $500$ and the number of test data set is $100$. Figure~\ref{fig-Ex1} shows the illustrated comparison between $y_t$ and $y_p$.
\begin{table}[t]
\caption{The exact value $y_t$, the predicted value of the ADM method, the predicted value of the FNN, the predicted value $y_p$ and the absolute error~(Error $|y_t - y_p|$) of the LDNN in several test points on $[0,1]$ domain for Experiment~1.}
\label{Ex1-table}
\begin{center}
\begin{tabular}{cccccc}
$x$&exact value&predicted value&predicted value&predicted value&Error\\
&($y_t=\mathrm{e}^{x}$)&by ADM&by FNN&($y_p$) by LDNN& \\
\\ \hline \\
$0.0$&$1.0$&$1.0002421$&$1.0006372$&$1.000000049$&$\expnumber{4.90000001}{-08}$\\
$0.2$&$1.22140276$&$1.2213538$&$1.2213246$&$1.221402765$&$\expnumber{4.99999997}{-09}$\\
$0.4$&$1.4918247$&$1.4919181$&$1.4919742$&$1.49182494$&$\expnumber{2.40000000}{-07}$\\
$0.6$&$1.8221188$&$1.8220339$&$1.8221628$&$1.82211831$&$\expnumber{4.90000000}{-07}$\\
$0.8$&$2.22554093$&$2.2255747$&$2.2255346$&$2.225540981$&$\expnumber{5.09999998}{-08}$\\
$1.0$&$2.71828183$&$2.717803$&$2.717317$&$2.71828179$&$\expnumber{4.00000002}{-08}$\\
\end{tabular}
\end{center}
\end{table}
\setlength\abovecaptionskip{-5pt}

\subsection{Experiment~2}
\begin{figure}[h]
\begin{center}
\begin{subfigure}[b]{0.3\textwidth}
         \begin{center}
         \includegraphics[width=\textwidth]{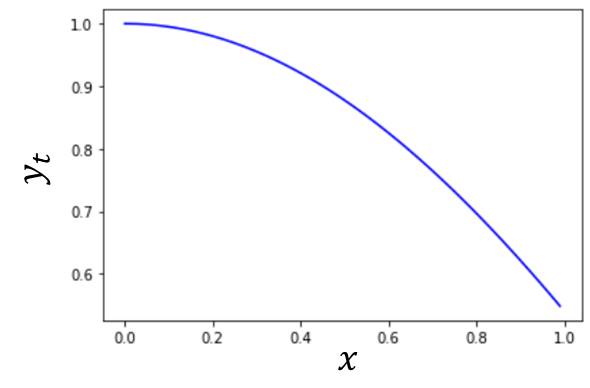}
         \label{fig:2-1}
         \end{center}
     \end{subfigure}
     \begin{subfigure}[b]{0.3\textwidth}
         \begin{center}
         \includegraphics[width=\textwidth]{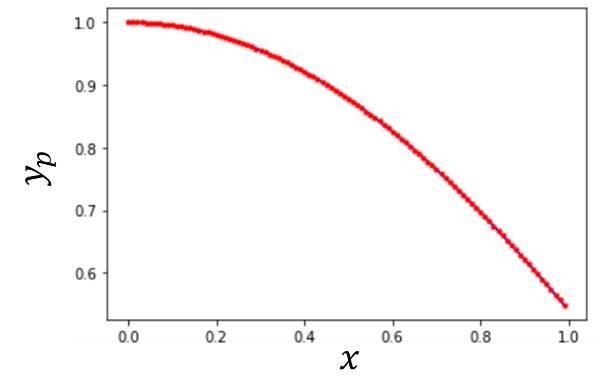}
         \label{fig:2-2}
         \end{center}
     \end{subfigure}
     \begin{subfigure}[b]{0.3\textwidth}
         \begin{center}
         \includegraphics[width=\textwidth]{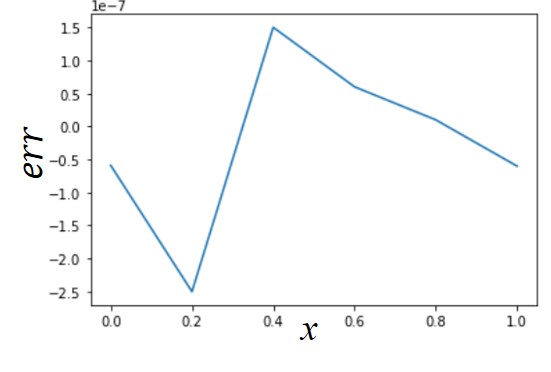}
         \label{fig:2-3}
         \end{center}
     \end{subfigure}
\end{center}
\caption{Results of Experiment~2. Exact solution $y_t(x)=\cos(x)$,  predicted solution $y_p(x)$ by LDNN and absolute error $err=|y_t(x)-y_p(x)|$.}
\label{fig-Ex2}
\end{figure}

Suppose that we have the following model \citep{razzaghi2002rationalized}:
\begin{equation}
\label{Ex2}
    y(x)=1 + \sin^2(x)+\int_{0}^{1} K(x,s)y^2(s) \mathrm{d}s, \quad x \in [0,1].
\end{equation}
where 
\begin{equation}
    K(x,s)=\big\{\begin{matrix}-3\sin(x-s),&0\leq s \leq x;\\
    0 & x\leq s \leq 1.\end{matrix}
\end{equation}
It has the exact solution $y(x)=\cos(x)$. The exact value, the predicted value and the absolute error~(Error) in several test points on $[0,1]$ domain are reported in Table~\ref{Ex2-table}. $50$ points of shifted Legendre quadrature points are applied for training LDNN. The number of train data set is $500$ and the number of test data set is $100$. Figure~\ref{fig-Ex2} shows the illustrated comparison between $y_t$ and $y_p$. 
\begin{table}[t]
\caption{The exact value $y_t$, the predicted value of the ADM method, the predicted value of the FNN, the predicted value $y_p$ and the absolute error~(Error $|y_t - y_p|$) of the LDNN in several test points on $[0,1]$ domain for Experiment~2.}
\label{Ex2-table}
\begin{center}
\begin{tabular}{cccccc}
\\
$x$&exact value&predicted value&predicted value&predicted value&Error\\
&($y_t=\cos(x)$)&by ADM&by FNN&($y_p$) by LDNN& \\
\\ \hline \\
$0.0$&$1.0$&$1.0003562$&$1.0006432$&$1.000000059$&$\expnumber{5.90000000}{-08}$\\
$0.2$&$0.98006658$&$0.97977763$&$0.97865423$&$0.98006683$&$\expnumber{2.50000000}{-07}$\\
$0.4$&$0.92106099$&$0.9210639$&$0.92105988$&$0.92106083$&$\expnumber{1.50000000}{-07}$\\
$0.6$&$0.82533561$&$0.82562345$&$0.82580132$&$0.82533555$&$\expnumber{6.00000000}{-08}$\\
$0.8$&$0.69670671$&$0.6963889$&$0.69647832$&$0.69670670$&$\expnumber{9.99999994}{-09}$\\
$1.0$&$0.54030231$&$0.5411298$&$0.54212091$&$0.54030237$&$\expnumber{6.00000001}{-08}$\\
\end{tabular}
\end{center}
\end{table}
\setlength\abovecaptionskip{-5pt}

\subsection{Experiment~3}
\begin{figure}[h]
\begin{center}
\begin{subfigure}[b]{0.3\textwidth}
         \begin{center}
         \includegraphics[width=\textwidth]{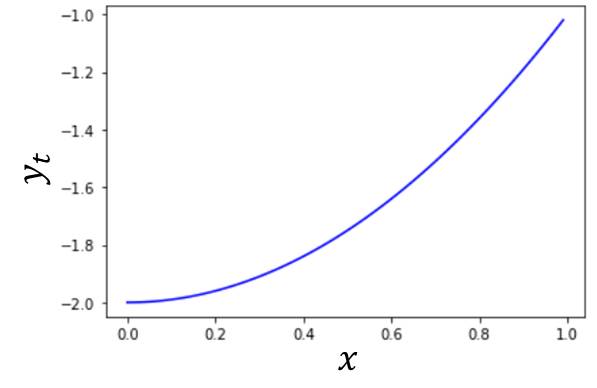}
         \label{fig:3-1}
         \end{center}
     \end{subfigure}
     \begin{subfigure}[b]{0.3\textwidth}
         \begin{center}
         \includegraphics[width=\textwidth]{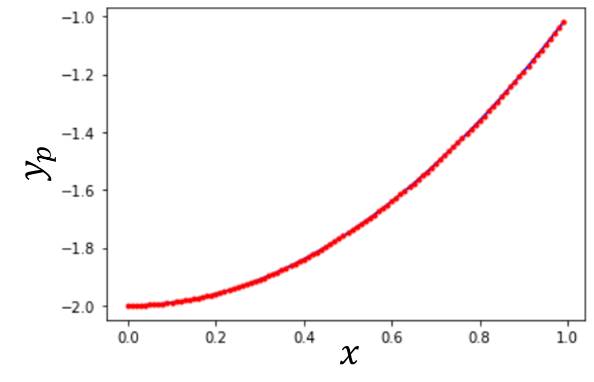}
         \label{fig:3-2}
         \end{center}
     \end{subfigure}
     \begin{subfigure}[b]{0.3\textwidth}
         \begin{center}
         \includegraphics[width=\textwidth]{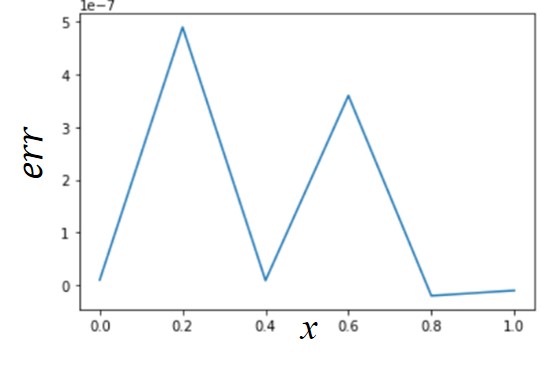}
         \label{fig:3-3}
         \end{center}
     \end{subfigure}
\end{center}
\caption{Results of Experiment~3. Exact solution $y_t(x)=x^2-2$,  predicted solution $y_p(x)$ by LDNN and absolute error $err=|y_t(x)-y_p(x)|$.}
\label{fig-Ex3}
\end{figure}

Suppose that we have the following model \cite{babolian2007chebyshev}:
\begin{equation}
\label{Ex3}
    y(x)=g(x) +\int_{0}^{x} (x-s)y^2(s) \mathrm{d}s + \int_{0}^{1} (x+s)y(s) \mathrm{d}s, \quad x \in [0,1].
\end{equation}
where 
\begin{equation}
    g(x)= - \frac{1}{30}x^6 + \frac{1}{3}x^4 - x^2 + \frac{5}{3}x - \frac{5}{4}
\end{equation}
It has the exact solution $y(x)=x^2 - 2$. Table~\ref{Ex3-table} illustrates the exact value, the predicted value and the absolute error~(Error) in several test points on $[0,1]$ domain. $50$ points of shifted Legendre quadrature points are applied for training LDNN. The number of train data set is $500$ and the number of test data set is $100$. Figure~\ref{fig-Ex3} represented the comparison between $y_t$ and $y_p$. 
\begin{table}[t]
\caption{The exact value $y_t$, the predicted value of the ADM method, the predicted value of the FNN, the predicted value $y_p$ and the absolute error~(Error $|y_t - y_p|$) of the LDNN in several test points on $[0,1]$ domain for Experiment~3.}
\label{Ex3-table}
\begin{center}
\begin{tabular}{llllll}
\\
$x$&exact value&predicted value&predicted value&predicted value&Error\\
&($y_t=x^2-2$)&by ADM&by FNN&($y_p$) by LDNN& \\
\\ \hline \\
$0.0$&$-2.0$&$-2.0001612$&$-2.0003422$&$-2.00000001$&$\expnumber{9.99999994}{-09}$\\
$0.2$&$-1.96$&$-1.960051$&$-1.9601312$&$-1.96000049$&$\expnumber{4.90000000}{-07}$\\
$0.4$&$-1.84$&$-1.8399543$&$-1.8396587$&$-1.840000009$&$\expnumber{8.99999986}{-09}$\\
$0.6$&$-1.64$&$-1.6400322$&$-1.640040$&$-1.64000036$&$\expnumber{3.60000000}{-07}$\\
$0.8$&$-1.36$&$-1.3599668$&$-1.35889879$&$-1.35999998$&$\expnumber{2.00000001}{-08}$\\
$1.0$&$-1.0$&$-0.9999476$&$-0.99987677$&$-0.99999999$&$\expnumber{1.00000001}{-08}$\\
\end{tabular}
\end{center}
\end{table}
\setlength\abovecaptionskip{-5pt}

\subsection{Experiment~4}
\begin{figure}[h]
\begin{center}
\begin{subfigure}[b]{0.3\textwidth}
         \begin{center}
         \includegraphics[width=\textwidth]{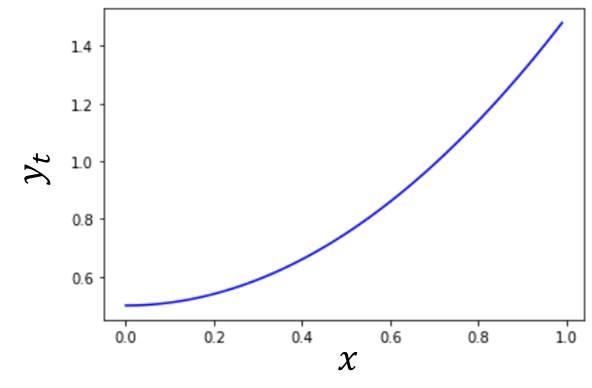}
         \label{fig:4-1}
         \end{center}
     \end{subfigure}
     \begin{subfigure}[b]{0.3\textwidth}
         \begin{center}
         \includegraphics[width=\textwidth]{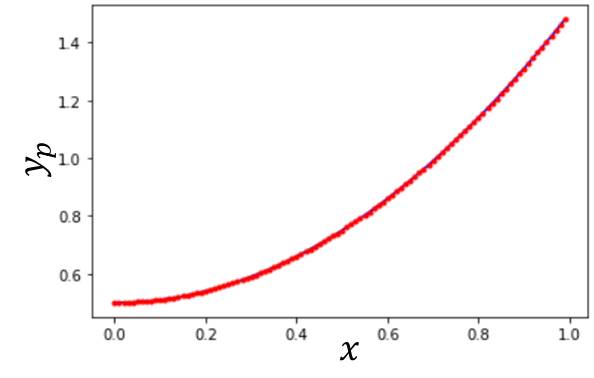}
         \label{fig:4-2}
         \end{center}
     \end{subfigure}
     \begin{subfigure}[b]{0.3\textwidth}
         \begin{center}
         \includegraphics[width=\textwidth]{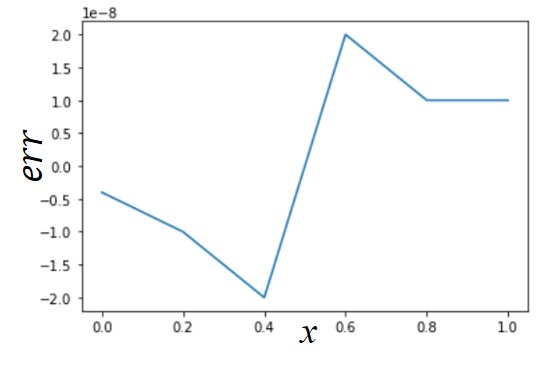}
         \label{fig:4-3}
         \end{center}
     \end{subfigure}
\end{center}
\caption{Results of Experiment~4. Exact solution $y_t(x)=x^2 + \frac{1}{2}$,  predicted solution $y_p(x)$ by LDNN and absolute error $err=|y_t(x)-y_p(x)|$.}
\label{fig-Ex4}
\end{figure}

Suppose that we have the following model \cite{hadizadeh2005numerical}:
\begin{equation}
\label{Ex4}
    y(x)=-\frac{1}{10}x^4 + \frac{5}{6}x^2 + \frac{3}{8} +\int_{0}^{x} \frac{1}{2x}y^2(s) \mathrm{d}s, \quad x \in [0,1].
\end{equation}
It has the exact solution $y(x)=x^2 + \frac{1}{2}$. The exact value, the predicted value and the absolute error~(Error) in several test points on $[0,1]$ domain are reported in Table~\ref{Ex4-table}. $50$ points of shifted Legendre quadrature points are applied for training LDNN. The number of train data set is $500$ and the number of test data set is $100$. Figure~\ref{fig-Ex4} shows the illustrated comparison between $y_t$ and $y_p$. 
\begin{table}[t]
\caption{The exact value $y_t$, the predicted value of the ADM method, the predicted value of the FNN, the predicted value $y_p$ and the absolute error~(Error $|y_t - y_p|$) of the LDNN in several test points on $[0,1]$ domain for Experiment~4.}
\label{Ex4-table}
\begin{center}
\begin{tabular}{llllll}
\\
$x$&exact value&predicted value&predicted value&predicted value&Error\\
&($y_t=x^2 + \frac{1}{2}$)&by ADM&by FNN&($y_p$) by LDNN& \\
\\ \hline \\
$0.0$&$0.5$&$0.50039285$&$0.50042379$&$0.50000004$&$\expnumber{4.00000000}{-09}$\\
$0.2$&$0.54$&$0.5400339$&$0.54006321$&$0.54000001$&$\expnumber{9.99999994}{-09}$\\
$0.4$&$0.66$&$0.6599865$&$0.6600465$&$0.66000002$&$\expnumber{2.00000000}{-08}$\\
$0.6$&$0.86$&$0.86001176$&$0.8598769$&$0.85999998$&$\expnumber{2.00000000}{-08}$\\
$0.8$&$1.14$&$1.1399317$&$1.13988365$&$1.13999999$&$\expnumber{9.99999994}{-09}$\\
$1.0$&$1.5$&$1.499708$&$1.4998377$&$1.49999999$&$\expnumber{9.99999994}{-09}$\\
\end{tabular}
\end{center}
\end{table}
\setlength\abovecaptionskip{-5pt}
\section{Conclusion}
Legendre deep neural network (LDNN) is introduced in this paper. LDNN and its application for solving nonlinear Volterra–Fredholm–Hammerstein integral equations (V-F-H-IEs) are proposed. LDNN includes two networks. The first network approximates the solution of a nonlinear V-F-H-IE $y(x)$ which has $ \mathcal{M}$-layers feed forward neural network structure. The first hidden layer of this has a orthogonal layer consists of Legendre polynomials as activation functions. The last network adjusts the output of the sooner network to fit to a desired equation form. The better performance of the network has been obtained by using Legendre Gaussian integration and automatic differentiation. Some experiments of nonlinear  V-F-H-IEs are given to investigate the reliability and validity of LDNN. The results show that this network is an efficient and has high accuracy.

\bibliography{Refs}

\end{document}